\documentclass[a4paper,11pt]{amsart}
\usepackage{amscd,amssymb,bbm}
\usepackage[all]{xy}

\newcommand{\gr}{\operatorname{gr}}

\newcommand{\Hom}{\operatorname{Hom}}

\newcommand{\Ker}{\operatorname{Ker}}

\newcommand{\Ext}{\operatorname{Ext}}

\newcommand{\op}{\operatorname{op}}

\newcommand{\fcog}{\operatorname{fcog}}
\newcommand{\lf}{\operatorname{l.f.}}

\newcommand{\G}{\Gamma}
\renewcommand{\L}{\Lambda}

\newcommand{\Z}{{\mathbb Z}}

\newcommand{\X}{{\mathcal X}}

\newcommand{\lan}{\langle}
\newcommand{\ran}{\rangle}

\newcommand{\Gr}{\operatorname{Gr}}

\newcommand{\K}{{\mathcal K}}

\newcommand{\pil}{\rightarrow}
\newcommand{\ainf}{A_{\infty}}
\newcommand{\extto}{\xrightarrow}

\newcommand{\Dgr}{\mathcal{D} \Gr}
\newcommand{\gsh}[1]{\lan #1 \ran}
\newcommand{\mfgr}{\operatorname{\mathfrak{gr}}}
\newcommand{\tmfgr}{\tilde{\operatorname{\mathfrak{gr}}}}
\newcommand{\md}[2]{#1 ^{( #2 )}}
\newtheorem{lem}{Lemma}[section]
\newtheorem{prop}[lem]{Proposition}
\newtheorem{cor}[lem]{Corollary}
\newtheorem{thm}[lem]{Theorem}
\theoremstyle{definition}
\newtheorem{defin}[lem]{Definition}
\newtheorem{example}[lem]{Example}
\newtheorem{remark}[lem]{Remark}

\begin{document}

\title{The Koszul dual of a weakly Koszul module}
\author{Dag Madsen}
\address{Institutt for matematiske fag, NTNU, NO--7491 Trondheim, Norway}
\email{dagma@math.ntnu.no}



\begin{abstract}
We study the so-called weakly Koszul modules and characterise their
Koszul duals. We show that the (adjusted) associated graded module
of a weakly Koszul module exactly determines the homology modules of
the Koszul dual. We give an example of a quasi-Koszul module which
is not weakly Koszul.
\end{abstract}

\maketitle

Let $\L=\bigoplus_{i \geq 0} \L_i$ be a Koszul algebra. Such
algebras have nice homological properties and appear in various
combinatorial and geometrical contexts. There is a corresponding
notion of Koszul modules (see section 1 below for definitions).
However, judging from the results in \cite{App}, especially Theorems
4.5 and 5.6 of that paper, for some important theoretical purposes
the class of Koszul modules is too small, and often the larger class
of so-called weakly Koszul modules is needed to complete the
picture. Let us briefly recall the definition, a more thorough
discussion will follow in section 1. Let $M$ be a finitely generated
graded $\L$-module. We say that $M$ is weakly Koszul if
$\Ext^{\ast}_\L (M,\L_0)=\bigoplus_{i \geq 0} \Ext^i_\L (M,\L_0)$ is
a Koszul module over the $\Ext$ algebra $\Ext^{\ast}_\L
(\L_0,\L_0)=\bigoplus_{i \geq 0} \Ext^i_\L (\L_0,\L_0)$. This is a
property satisfied by the Koszul modules themselves, so Koszul
modules are weakly Koszul.

In the present paper we investigate the behavior of weakly Koszul
modules under the Koszul duality functor of \cite{Bei}, which is a
functor on the level of derived categories. We show that Koszul
duals of weakly Koszul modules can be characterised in terms of
their homology (Theorem \ref{char}). We also show that the Koszul
duals of two weakly Koszul modules have isomorphic homology modules
if and only if the two weakly Koszul modules have isomorphic
(adjusted) associated graded modules with respect to the radical
filtration (Corollary \ref{sist}). When investigating weakly Koszul
modules it is therefore relevant to study the objects in the derived
category with given homology modules. The language of
$\ainf$-modules \cite{Kel} is particularly well suited for this
purpose, and we exploit this in the last section.

Our description of Koszul duals of weakly Koszul modules can be used
as a basis for studying more general classes of modules. We say that
$M$ is quasi-Koszul if $\Ext^{\ast}_\L (M,\L_0)$ is generated in
degree $0$. Weakly Koszul modules are quasi-Koszul, and quasi-Koszul
modules generated in a single degree are weakly Koszul. Based on
this evidence, one could be led to believe that the two notions are
equivalent. Such speculations can now be laid to rest, as we present
a counterexample (Example \ref{quasi}). Generalising in another
direction we consider modules $M$ such that $\Ext^{\ast}_\L
(M,\L_0)$ is weakly Koszul. These modules have a surprising
property, compare Theorem \ref{eweak} and Example \ref{weak}.

The contents of the different sections are as follows. In section 1
we give the basic definitions of Koszul algebras and modules, weakly
Koszul modules and other related classes of modules. We also recall
the fundamental dualities and equivalences present in this setting.
Section 2 explains the concept of the "adjusted" associated graded
module of $M$. In section 3 we give several results and formulas
concerning the homology of the Koszul dual of a weakly Koszul
module. In the two last sections we discuss how to find the object
itself, not just its homology. Section 4 gives a method for
straightforward computations, and we use this to produce some
important (counter)examples. In section 5 we take a more systematic
approach. We show how to use $\ainf$-module structures to classify
all objects with given homology modules. Via Koszul duality, this
classifies all weakly Koszul modules with a given adjusted
associated graded module.

\section{dualities and equivalences}

Let $k$ be a field and $\L=\bigoplus_{i \geq 0}\L_i$ be a graded
$k$-algebra. We assume $\dim_k \L_i < \infty$ for all $i \geq 0$ and
$\L_0 \simeq k \times \ldots \times k$ as rings. We denote by $J$
the graded Jacobson radical $J=\bigoplus_{i \geq 1}\L_i$. We denote
by $\Gr \L$ the category of graded $\L$-modules $M=\bigoplus_{i \in
\Z} M_i$ with degree $0$ morphisms. By $\lf \L$ we denote the full
subcategory of locally finite modules, that is modules with $\dim_k
M_i < \infty$ for all $i \in \Z$. Important (full) subcategories of
$\lf \L$ are the category of finitely generated modules $\gr \L$ and
the category of finitely cogenerated modules $\fcog \L$.

Let $M$ be a graded $\L$-module. Its \emph{graded dual} $DM$ is
defined to be the graded $\L^{\op}$-module with graded parts
$(DM)_{-i}=\Hom_k(M_i,k)$ and (graded parts of) module structure
maps $k$-dual to those of $M$. Also using $k$-dual maps on morphisms
we can make $D$ into a contravariant functor $D \colon \Gr \L \pil
\Gr \L^{\op}$. When restricted to locally finite $\L$-modules, the
functor $D$ is a duality $D \colon \lf \L \pil \lf \L^{\op}$.

The following lemma, which can for instance be found in
\cite[2.4.7]{Grad}, gives a useful connection between graded and
ungraded Ext groups.

\begin{lem}\label{gung}
Suppose $M$ is a finitely generated graded $\L$-module which has a
projective resolution such that all syzygies are finitely generated.
Let $L$ be any graded $\L$-module.

Then for every $j \geq 0$ we have an isomorphism
$$\Ext^j_{\L}(M,N) \simeq \bigoplus_{i \in \Z} \Ext^j _{\Gr \L}(M,N
\gsh i)$$ functorial in $M$ and $N$.
\end{lem}

Typical examples of modules $M$ which satisfy the hypothesis of the
lemma are finitely generated graded $\L$-modules over a Noetherian
algebra $\L$.

We say that a graded $\L$-module $M$ is a \emph{Koszul module} if
$M$ is finitely generated and $\Ext^j_{\Gr \L}(M,\L_0 \gsh i) \neq
0$ implies $i=j$. In particular we require that $\Hom_{\Gr
\L}(M,\L_0 \gsh i) \neq 0$ implies $i=0$, so $M$ must be generated
in degree $0$.

\begin{remark}
If $L$ is a graded $\L$-module that is locally finite and bounded
below (for instance if $L$ is finitely generated), then there is a
projective cover $P \extto{^f} L \pil 0$, where $P$ is a projective
module which is locally finite and bounded below. It follows that
$\Ker f$ is also locally finite and bounded below. Such a module $L$
will therefore have a minimal graded projective resolution. The
$\Ext$ condition in the definition of a Koszul module is a condition
on the degrees where the projective modules in the resolution are
generated. More precisely  a finitely generated module $M$ is Koszul
if and only if $M$ has a graded projective resolution
$$\ldots \pil P_n \pil \ldots \pil P_2 \pil P_1 \pil P_0 \pil M \pil 0$$
with $P_i$ finitely generated in degree $i$ for all $i \geq 0$.
\end{remark}

A locally finite $\L$-module $M$ is called a \emph{coKoszul module}
if $DM$ is a Koszul $\L^{\op}$-module. We denote the full
subcategory of Koszul $\L$-modules by $\K(\L)$ and the full
subcategory of coKoszul $\L$-module by $c\K(\L)$. The functor $D$
restricts to dualities on subcategories in the way shown by the
following diagram.

$$\xymatrix{\lf \L \ar[r]^D  & \lf \L^{\op}\\
\gr \L \ar[r]^D \ar@{^{(}->}[u] & \fcog \L^{\op} \ar@{^{(}->}[u]\\
\K(\L) \ar[r]^D \ar@{^{(}->}[u] & c\K(\L^{\op}) \ar@{^{(}->}[u]
}$$

The \emph{$i$th graded shift} of $M$, denoted $M \gsh i$, is the
module with graded parts $(M \gsh i)_n=M_{n-i}$ and module structure
inherited from $M$. If $M$ is a module generated in a single degree
$i$, we define $\overline M=M \gsh {-i}$. So $\overline M$ is
generated in degree $0$. If $\overline M$ is Koszul, we say that $M$
has a \emph{linear resolution}.

The algebra $\L$ is called a \emph{Koszul algebra} if $\L_0$ is a
Koszul $\L$-module. One can prove that $\L$ is a Koszul algebra if
and only if $\L^{\op}$ is a Koszul algebra. Suppose now and for the
rest of the paper that $\L$ is a Koszul algebra. Let
$\G=\bigoplus_{i \geq 0}\Ext^i_\L(\L_0,\L_0) \simeq \bigoplus_{i
\geq 0}\Ext^i_{\Gr \L} (\L_0,\L_0 \gsh i)$. A fundamental theorem
\cite[1.2.5]{Bei} states that $\G$ is also a Koszul algebra. The
algebra $\G$ is called the \emph{Koszul dual} of $\L$. The Koszul
dual of $\G$ is isomorphic to $\L$ as graded algebras.

Another fundamental theorem states that there is an equivalence
between certain triangulated subcategories of the corresponding
(unbounded) derived categories. The category $\mathcal
D^{\downarrow} (\L^{\op})$ can be viewed as the full subcategory of
$\Dgr \L^{\op}$ formed by all objects $M$ with the property that
$(H^i M)_j=0$ when $i \ll 0$ or $i+j \gg 0$. Similarly, the category
$\mathcal D^{\uparrow} (\G)$ is the full subcategory of $\Dgr \G$
formed by the objects $N$ with the property that $(H^i N)_j=0$ when
$i \gg 0$ or $i+j \ll 0$. The theorem \cite[2.12.1]{Bei} states that
there is an equivalence of triangulated categories $G \colon
\mathcal D^{\downarrow} (\L^{\op}) \pil \mathcal D^{\uparrow} (\G)$.
A concrete description of the functor $G$ is given in \cite{Bei},
and we shall use this description for computations in section
\ref{diff}. (To avoid confusion we point out that the functor is
called $K$ in \cite{Bei}.)

Another description of the functor $G$ follows from the theory of
"lifts" in \cite[7.3,10]{Dgc}. There is a bigraded
$\L^{\op}$-$\G$-bimodule complex $X$, in degree $(\ast, i)$
quasi-isomorphic to $\L_0^{\op} \gsh {-i} [-i]$ as a complex of left
$\L^{\op}$-modules, such that the functor $R\Hom(X,-) \colon \Dgr
\L^{\op} \pil \Dgr \G$ when restricted to $D^{\downarrow}
(\L^{\op})$ is isomorphic to $G$. In this paper we will not attempt
to describe the bimodule $X$ further. For a discussion of to which
extent the category equivalence determines $X$ , we refer to
\cite[7]{Dgc}.

In the usual way we view modules as stalk complexes concentrated in
degree $0$. The category $\mathcal D^{\downarrow} (\L^{\op})$
contains all finitely cogenerated modules. From the isomorphism $G
\simeq R\Hom(X,-)$ we get the following result.

\begin{prop}\label{rhom}
Let $M$ be a finitely cogenerated $\L^{\op}$-module. Then
\begin{itemize}
\item[(a)] $(H^j G(M))_i \simeq \Ext^{i+j}_{\Gr \L^{\op}}(\L^{\op}_0,M \gsh i).$
\item[(b)] $G(M \gsh i) \simeq (G(M)) \gsh {-i} [-i]$.
\end{itemize}
\end{prop}

It follows from part (a) that if $M$ is a coKoszul $\L$-module, then
$G(M)=\bigoplus_{i \geq 0}\Ext^i_{\Gr \L^{\op}} (\L_0,M \gsh i)$ is
a $\G$-module. It is possible to show \cite[5.1]{Gre} that in this
case $G(M)$ is a Koszul $\G$-module. From Lemma \ref{gung} we get an
isomorphism of $\G$-modules $\bigoplus_{i \geq 0}\Ext^i_{\Gr
\L^{\op}} (\L_0,M \gsh i) \simeq \bigoplus_{i \geq
0}\Ext^i_{\L^{\op}} (\L_0,M)$ (functorial in $M$). So when $G$ is
restricted to $c\K (\L^{\op})$ it is isomorphic to the functor
$E=\bigoplus_{i \geq 0} \Ext^i_{\L^{\op}}(\L_0,-)$. The relation
between the various categories and functors is summed up in the
following diagram.

$$\xymatrix{& \Dgr \L^{\op} \ar[rr]^{R \Hom(X,-)} && \Dgr \G \\
& \mathcal D^{\downarrow} (\L^{\op}) \ar@{^{(}->}[u]
\ar[rr]^{\sim}_G &&
\mathcal D^{\uparrow} (\G) \ar@{^{(}->}[u]\\
\K (\L) \ar[r]^D & c\K (\L^{\op}) \ar@{^{(}->}[u] \ar[rr]^{\sim}_E
&& \K (\G) \ar@{^{(}->}[u] }$$

Let $\check E=ED$. Then $\check E \colon \K(\L) \pil \K(\G)$ is also
a duality and can be described as $\check E=\bigoplus_{i \geq 0}
\Ext^i_\L(-,\L_0)$. This functor can be applied to any graded
module, so we view $\check E$ with this description as a functor
$\check E \colon \Gr \L \pil \Gr \G$. Note that this functor forgets
the $\L$-grading, so $\check E(M) \simeq \check E(M \gsh i)$ for any
graded $\L$-module $M$ and $i \in \Z$. We have the following lemma
concerning local finiteness of $\check E(M)$.

\begin{lem}\label{hatem}
Let $M$ be a finitely generated graded $\L$-module.
\begin{itemize}
\item[(a)] Suppose $GD(M)$ has bounded homology and that $\dim_k (H^i
GD(M))_j < \infty$ for all $i,j \in \Z$. Then $M$ has a projective
resolution such that all syzygies are finitely generated.
\item[(b)] The $\G$-module $\check E(M)$ is locally finite if and only if $M$ has
a projective resolution such that all syzygies are finitely
generated. Moreover, in this case $\check E(M) \simeq \bigoplus_{i
\in \Z} (H^i GD(M)) \gsh i$ as graded $\G$-modules.
\end{itemize}
\end{lem}

\begin{proof}
Suppose $M$ has a minimal graded projective resolution
$$\ldots \pil P_n \pil \ldots \pil P_2 \pil P_1 \pil P_0 \pil M \pil 0.$$

(a): From minimality it follows that $\Hom_{\Gr \L}(P_j,\L_0 \gsh l)
\simeq \Ext^j_{\Gr \L}(M,\L_0 \gsh l) \simeq \Ext^j_{\Gr
\L^{\op}}(\L_0^{\op},DM \gsh l) \simeq (H^{l-j} GD(M))_l$ for all $l
\in \Z$, $j \geq 0$. For any given $j \geq 0$, since $GD(M)$ has
bounded homology, this is nonzero only for a finite number of values
of $l$ (and finite dimensional in those cases). Therefore $P_j$ is
finitely generated.

(b): From minimality of the above sequence it follows that $\check
E(M)_j=\Ext^j_\L(M,\L_0) \simeq \Hom_\L(P_j,\L_0)$ for all $j \in
\Z$.

Suppose $P_j$ is not finitely generated for some $j \geq 0$. Since
$P_j$ as a graded module is locally finite and bounded below, we
must have $\dim_k (\Hom_{\Gr\L}(P_j,\L_0))=\infty$. Since
$\Hom_{\Gr\L}(P_j,\L_0) \subseteq \Hom_\L(P_j,\L_0)$, the module
$\check E(M)$ is not locally finite.

Suppose $P_j$ is finitely generated for all $j \geq 0$. Then $\check
E(M)_j \simeq \Hom_\L(P_j,\L_0)$ is finite dimensional for all $j
\geq 0$ (and zero for $j<0$), so $\check E(M)$ is locally finite. By
Lemma \ref{gung} we have $\check E(M)_j \simeq \bigoplus_{i \in \Z}
\Ext^j_{\Gr \L}(M,\L_0 \gsh i)$ for all $j \in \Z$. From Proposition
\ref{rhom} we get $(\bigoplus_{i \in \Z} (H^i GD(M)) \gsh
i)_j=\bigoplus_{i \in \Z} ((H^i GD(M)) \gsh {j-i}) \simeq
\bigoplus_{i \in \Z} \Ext^j_{\Gr \L^{\op}}(\L_0^{\op},DM \gsh i)
\simeq \check E(M)_j$ for all $j \in \Z$. Also the module structure
is preserved, so $\check E (M) \simeq \bigoplus_{i \in \Z} (H^i
GD(M)) \gsh i$ as graded $\G$-modules.
\end{proof}

We are now ready to define weakly Koszul modules.

\begin{defin}
A finitely generated graded $\L$-module $M$ is called a \emph{weakly
Koszul module} if $\check E(M)$ is a Koszul $\G$-module.

A slightly weaker condition is that $\check E(M)$ is generated in
degree $0$. A finitely generated graded $\L$-module $M$ satisfying
this condition is called a \emph{quasi-Koszul module}.
\end{defin}

\begin{remark}
There are other equivalent ways of defining weakly Koszul modules.
For instance, a finitely generated graded $\L$-module
$M=\bigoplus_{i \in \Z} M_i$ is weakly Koszul if and only if for all
$j \in \Z$, the submodule of $M$ generated by $M_j$ has a linear
resolution. Quasi-Koszul and weakly Koszul (originally under the
name \emph{strongly quasi-Koszul}) modules were introduced in
\cite{Coc}. In addition to the mentioned paper \cite{App}, weakly
Koszul modules are also studied in the paper \cite{Nwk}.
\end{remark}

Suppose $M$ is generated in a single degree. In this case it can be
shown that $M$ is quasi-Koszul if and only if $M$ is weakly Koszul,
which is again equivalent to $M$ having a linear resolution. At the
end of section \ref{diff} we give an example of a module $M$
generated in multiple degrees which is quasi-Koszul but not weakly
Koszul.

We denote the category of weakly Koszul $\L$-modules by $w\K (\L)$.
It is closed under direct summands and finite direct sums. We call a
module dual under $D$ to a weakly Koszul module a \emph{weakly
coKoszul} module. The corresponding category of weakly coKoszul
$\L^{\op}$-modules we denote by $wc\K(\L^{\op})$.

The essential image of $wc\K(\L^{\op})$ under $G \simeq R\Hom(X,-)$
we denote by $\X$. Our aim is to describe this category. It is a
full subcategory of $D^{\uparrow} (\G)$ as we see from the following
diagram.

$$\xymatrix{& \Dgr \L^{\op} \ar[rr]^{R\Hom(X,-)}  && \Dgr \G \\
& \mathcal D^{\downarrow} (\L^{\op}) \ar@{^{(}->}[u] \ar[rr]^{\sim}
&&
\mathcal D^{\uparrow} (\G) \ar@{^{(}->}[u]\\
w\K (\L) \ar[r]^D & wc\K(\L^{\op}) \ar@{^{(}->}[u] \ar[rr]^{\sim}
&& \X \ar@{^{(}->}[u]\\
\K (\L) \ar[r]^D & c\K (\L^{\op}) \ar@{^{(}->}[u] \ar[rr]^{\sim}_E
&& \K (\G) \ar@{^{(}->}[u] }$$

\section{associated graded module}
In this section we explain some technicalities concerning
filtrations of finitely generated modules and the associated graded
modules.

If $M$ is a finitely generated $\L$-module, then its associated
graded module (with respect to the radical filtration) is $$\mfgr
(M)=\bigoplus_{i \geq 0}J^i M/J^{i+1} M.$$ This is also a finitely
generated graded $\L$-module. If $M$ is generated in a single
degree, then $\mfgr (M) \simeq \overline M$.

Suppose $M$ is finitely generated in degrees $j_0 < \ldots < j_p$.
We always assume that the set of generators is minimal, in other
words $M$ is finitely generated and $\Hom(M,\L_0 \gsh i) \neq 0$ if
and only if $i \in \{j_0, \ldots, j_p \}$. In \cite{App} we find the
following result.

\begin{prop}\cite[2.5]{App}
Suppose $M$ is a finitely generated module generated in degrees $j_0
< \ldots < j_p$. Let $K^{(0)}$ be the submodule of $M$ generated by
$M_{j_0}$. Then there is a split-exact sequence of $\L$-modules
$$0 \pil \mfgr(K^{(0)}) \pil \mfgr(M) \pil \mfgr(M/K^{(0)}) \pil  0.$$
\end{prop}

We define $\md M 0=M $, and in general for all $0 < i \leq p$ we
define $\md M i = \md M {i-1} / \md K {i-1}$ and let $\md K i$
denote the the submodule of $\md M i$ generated in degree $j_i$ (the
"highest degree"). For each $i$ the module $\md M i$ is generated in
degrees $j_i< \ldots <j_p$. In particular we have $\md M p=\md K p$.

With this notation we get the following corollary.

\begin{cor}\label{gr}
If $M$ is a finitely generated module generated in degrees $j_0 <
\ldots < j_p$, then
$$\mfgr(M) \simeq \mfgr (\bigoplus_{i=0}^p \md K i) \simeq
\bigoplus_{i=0}^p  \overline {\md K i}.$$
\end{cor}

Mention should now be made of the following theorem in \cite{App}.

\begin{thm}\cite[2.5]{App}
A finitely generated module $M$ is a weakly Koszul module if and
only if $\mfgr (M)$ is Koszul.
\end{thm}

Motivated by Corollary \ref{gr} we define the following "adjusted"
version of the associated graded module. If $M$ is a finitely
generated module, then $\tmfgr (M)$ is defined to be the module
$$\tmfgr (M) \simeq \bigoplus_{i=0}^p \md K i.$$

With this definition $M$ and $\tmfgr (M)$ are generated in the same
degrees, but each indecomposable summands of $\tmfgr (M)$ is
generated in a single degree. We also have $$\tmfgr (\md M s) \simeq
\bigoplus_{i=s}^p \md K i$$ whenever $0 \leq s \leq p$.

\section{homology of the Koszul dual}

We now return to the question of describing the objects in $\X$. The
following theorem shows that such objects can be characterised by
their homology.

\begin{thm}\label{char}
Let $N$ be an object in $\Dgr \G$.

Then $N \in \X$ if and only if
\begin{itemize}
\item[(i)] $N$ has bounded homology and
\item[(ii)] for all $i \in \Z$, $H^i N$ is generated in degree $-i$ and has
a linear resolution.
\end{itemize}
\end{thm}

\begin{proof}
If $N \in \X$, then $N \simeq GD(M)$ for a weakly Koszul $\L$-module
$M$. Also $\check E(M) \simeq \bigoplus_{i \geq 0} (H^i N) \gsh i$
is a Koszul $\G$-module and in particular it is generated in degree
$0$. Therefore for each $i \in \Z$, we have that  $(H^i N) \gsh i$
is Koszul and $(H^i N)$ is generated in degree $-i$. Since $\check
E(M)_0$ is finite dimensional over $k$, it means that $N$ has
bounded homology.

Let $N$ be an object in $\Dgr \G$ with bounded homology and suppose
$H^i N$ is generated in degree $-i$ for all $i \in \Z$. Choose a
representing complex for $N$ such that $N^i=0$ for $i \leq a$ and
$i>b$, for suitable integers $a,b$ with $a < b$. For each integer
$p$, the \emph{soft truncation} $\tau_{\leq p}N$ of $N$ at $p$ is
defined by
$$(\tau_{\leq p}N)^i= \left \{
\begin{array}{c c}
N^i & \text { if } i < p\\
\Ker d^i & \text { if } i=p\\
0 & \text { if } i > p
\end{array}
\right .$$ Its homology is given by
$$H^i(\tau_{\leq p}N)= \left \{
\begin{array}{c c}
H^i N & \text { if } i \leq p\\
0 & \text { if } i > p
\end{array}
\right .$$ We have a filtration $0=\tau_{\leq a} N \subseteq \ldots
\subseteq \tau_{\leq b} N=N$. All these objects are in $D^{\uparrow}
(\G)$. Consider the triangles in $D^{\uparrow} (\G)$
$$\tau_{\leq i-1} N \pil \tau_{\leq {i}} N \pil Y_i  \pil (\tau_{\leq i-1}
N) [1]$$ for all $i$ with $a<i \leq b$. Here $Y_i$ has nonzero
homology only possibly in degree $i$ and $H^i(Y_i) \simeq H^i N$. By
assumption, $H^i N \simeq L_i \gsh {-i} $ for some Koszul
$\G$-module $L_i$. So there is a Koszul $\L$-module $K_i$ such that
$GD(K \gsh {-i}) \simeq G ((DK) \gsh i) \simeq L_i \gsh {-i} [-i]
\simeq Y_i$.

So for all $i$ we have that $Y_i$ is isomorphic to $G$ of a finitely
cogenerated module (viewed as a stalk complex). Let $F \colon
D^{\uparrow} (\G) \pil D^{\downarrow} (\L^{\op})$ denote a
quasi-inverse of $G$. By induction (starting with $\tau_{\leq a}
N=0$), using the triangles above, we get for all $a < i \leq b$ that
$F (\tau_{\leq {i}} N)$ is a module and there are exact sequences
$$0 \pil F (\tau_{\leq {i-1}} N) \pil F (\tau_{\leq {i}} N) \pil
F(Y_i) \pil 0.$$ The modules $F(Y_i)$ are finitely cogenerated and
again by induction every $F (\tau_{\leq {i}} N)$ is finitely
cogenerated. So in particular $N=\tau_{\leq b} N$ is isomorphic to
$GD(M)$ for the finitely generated $\L$-module $M=DF(N)$. From Lemma
\ref{hatem} we get $\check E (M) \simeq \bigoplus_{i \in \Z} (H^i
GD(M)) \gsh i$ and by assumption this is a Koszul $\G$-module. But
then $M$ is a weakly Koszul module by definition.
\end{proof}

\begin{remark}
In Example \ref{quasi} we give an example of a graded $\L$-module
$M$ with the property that $N=GD(M)$ has bounded homology and $H^i
N$ is generated in degree $-i$ for all $i \in \Z$, but there is an
$i \in Z$ such that $H^i N$ does not have a linear resolution. This
means that the module $M$ is quasi-Koszul but not weakly Koszul.
\end{remark}

Let $M$ be a weakly Koszul module. We next try to find formulas for
the homology of $GD(M)$. We start with the simple case when the
module is a graded shift of a Koszul module.

\begin{prop}\label{lin}
Let K be a Koszul module. Then

$$H^{-n} GD (K \gsh i) \simeq \left \{
\begin{array}{c c} \check E(K) \gsh i  & \text{ if } n=i
\\ 0 & \text{ if } n \neq i
\end{array}
\right .
$$
\end{prop}

\begin{proof}
We have $H^{-n} GD (K \gsh i) \simeq H^{-n} G ((DK) \gsh {-i})
\simeq H^{-n} (GDK) [i] \gsh i \simeq H^{-n+i}(GDK) \gsh i \simeq
H^{-n+i} (\check E(K)) \gsh i$.
\end{proof}

The following proposition from \cite{App} will help us resolve the
general case.

\begin{prop}\cite[2.4]{App}\label{induc}
Let $M$ be a weakly Koszul module generated in degrees $j_0 < \ldots
< j_p$. Let $K^{(0)}$ be the submodule of $M$ generated by
$M_{j_0}$.

Then $\md K 0$ has a linear resolution and $\md M 1=M/\md K 0$ is
weakly Koszul.
\end{prop}

Keeping the notation from the previous section we have the following
obvious corollary.

\begin{cor}
Let $M$ be a weakly Koszul module generated in degrees $j_0 <
\ldots <j_p$.

Then for each $0 \leq i \leq p$, the module $\md K i$ has a linear
resolution and $\md M i$ is weakly Koszul.
\end{cor}

We know that the homology of $GD(M)$ is bounded and in each degree
it is of the same form as in Proposition \ref{lin}. Therefore there
must exist a module $\tilde M$, being a finite direct sum of modules
with linear resolutions, such that $H^n GD M \simeq H^n GD (\tilde
M)$ for all $n$. But what is this module  $\tilde M$? The next
proposition shows that the answer is the (adjusted) associated
graded module of $M$.

\begin{prop}\label{tmfgr}
Let $M$ be a weakly Koszul module. Then $$H^n GD M \simeq H^n GD (\tmfgr (M))$$
for all $n \in \Z$.
\end{prop}

\begin{proof}
We prove $H^n GD \md M s \simeq H^n GD (\tmfgr (\md M s))$ for all
$0 \leq s \leq p$ by downward induction on $s$ going from $\md M p=
\md K p$ to $\md M 0=M$.

Since $\md M p$ is generated in a single degree, we have $\md M p
\simeq \tmfgr (\md M p)$, so this case is clear.

From each exact sequence
$$0 \pil \md K s \pil \md M s \pil
\md M {s+1} \pil 0$$ with $0 \leq s < p$ we get a triangle
$$ D \md M {s+1} \pil D \md M s
\pil D \md K s \pil D \md M {s+1} [1]$$ in $\Dgr \L^{\op}$. Applying
$G$ to this triangle we get a triangle
$$ GD \md M {s+1} \pil GD \md M s
\pil GD \md K s \pil GD \md M {s+1} [1]$$ in $\Dgr \G$. We have a
long-exact sequence in homology
\begin{multline*}
\ldots \pil H^{n-1} GD \md K s \pil H^n GD \md M {s+1} \pil H^n GD
\md M s\\ \pil H^n GD \md K s \pil H^{n+1} GD \md M {s+1}  \pil
\ldots
\end{multline*}

We assume that $H^n GD \md M {s+1} \simeq H^n GD (\tmfgr (\md M
{s+1}))$ for a given $s$. The module $\md M {s+1}$ is generated in
degrees $j_{s+1} < \ldots <j_p$ and the same is true for $\tmfgr
(\md M {s+1}) \simeq \bigoplus_{i={s+1}}^p \md K i$. Since the $\md
K i$ have linear resolutions, Proposition \ref{lin} says that $H^n
GD \md M {s+1}$ is nonzero only for $n \in \{-j_p, \ldots, -j_{s+1}
\}$. In particular $H^n GD \md M {s+1} =0$ when $n \geq -j_s$. Also
by Proposition \ref{lin} $H^n GD \md K s \neq 0$ if and only if
$n=-j_s$.  Using these facts and the isomorphism $$\tmfgr (\md M s)
\simeq \md K s \oplus \tmfgr (\md M {s+1})$$ we get
$$H^n GD \md M s \simeq H^n GD \md K s \simeq H^n GD (\tmfgr (\md M s))$$
when $n \geq -j_s$ and
$$H^n GD \md M s \simeq H^n GD \md M {s+1} \simeq H^n GD (\tmfgr (\md M
{s+1})) \simeq H^n GD (\tmfgr (\md M s))$$ when $n<-j_s$. This
finishes the induction step.
\end{proof}

As a corollary we have the the following.

\begin{cor}\label{sist}
Let $M$ and $M'$ be two weakly Koszul modules. Then $$\tmfgr (M)
\simeq \tmfgr (M')$$ if and only if $$H^n GD M \simeq H^n GD M'$$
for all $n \in \Z$.
\end{cor}

\begin{proof}
The modules $\tmfgr (M)$ and $\tmfgr (M')$ are both direct sums of
modules with linear resolutions. From Proposition \ref{lin} it
follows that $\tmfgr (M) \simeq \tmfgr (M')$ if and only if $H^n GD
(\tmfgr (M)) \simeq H^n GD (\tmfgr (M'))$ for all $n \in \Z$. But by
Proposition \ref{tmfgr} this is true if and only if $H^n GD M \simeq
H^n GD M'$ for all $n \in \Z$.
\end{proof}

Combining Proposition \ref{tmfgr} with Proposition \ref{lin} we get
the following formula. Here $j^{-1}n$ is the number $i$ such that
$j_i=n$.

\begin{cor}\label{form}
Let $M$ be a weakly Koszul module generated in degrees $J=\{j_0,
\ldots, j_p\}$. Then $$H^{-n} GD (M) \simeq \left \{
\begin{array}{c c} \check E (\md K {j^{-1}n}) \gsh n  & \text{ if } n \in J
\\ 0 & \text{ if } n \notin J
\end{array}
\right .
$$
\end{cor}

In \cite{Ele}, the authors ask which finitely generated graded
$\L$-modules $M$ have the property that $\check E(M)$ is weakly
Koszul. We present the following proposition as a first step towards
understanding such modules.

\begin{thm}\label{eweak}
Suppose $M$ is a finitely generated graded $\L$-module such that
$\check E(M)$ is weakly Koszul. Then $\check E(M)$ has a direct
summand which is Koszul.
\end{thm}

\begin{proof}
Without loss of generality we may assume that $M$ is generated in
degrees $0=j_0 < j_1 < \ldots < j_p$. In this case $H^n GD(M)=0$ for
$n>0$ and there is a triangle $$\tau_{\leq -1}GD(M) \pil GD(M) \pil
H^0 GD(M) \pil \tau_{\leq -1}GD(M)[1]$$ in $D^{\uparrow} (\G)$.
Since $H^0 GD(M)$ is a direct summand of $\check E(M)$ by
\ref{hatem}, it is weakly Koszul by assumption. Since $(H^0 GD(M))_i
\simeq \Ext^i_{\Gr \L^{\op}}(\L^{\op}_0,DM \gsh i) \simeq
\Ext^i_{\Gr \L}(M, \L_0 \gsh i)$, we know that $(H^0 GD(M))$ has
support only in non-negative degrees. Also $(H^0 GD(M))_0 \simeq
\Hom_{\Gr \L}(M,\L_0) \neq 0$.

Let $F \colon D^{\uparrow} (\G) \pil D^{\downarrow} (\L^{\op})$
denote a quasi-inverse of $G$. If $K$ is a Koszul $\G$-module, then
$F(K)$ is a coKoszul $\L^{\op}$-module. Let $N$ be a weakly Koszul
$\G$-module generated in degrees $J=\{h_0, \ldots, h_p\}$. Since
$F(K \gsh i) \simeq (F(K) \gsh {-i}) [-i]$ for all $i \in \Z$, by
induction using Proposition \ref{induc} we get that $H^i(F(N))$ is
cogenerated in degree $-i$ for $i \in \Z$. Also $H^i(F(N)) \neq 0$
if and only if $i \in J$.

Let $N=H^0 GD(M)$. Applying $F$ to the triangle above, we get the
triangle $$F(\tau_{\leq -1}GD(M)) \pil DM \pil F(N) \pil
F(\tau_{\leq -1}GD(M))[1]$$ in $D^{\downarrow} (\L^{\op})$. Since
$N$ is generated in degrees $0=h_0 < h_1 < \ldots <h_p$, we have
$H^i F(N)=0$ for $i<0$. Since $DM$ is concentrated in homological
degree $0$, there is an isomorphism $H^i F(N) \simeq
H^{i+1}(F(\tau_{\leq -1}GD(M)))$ when $i>0$. Now $\tau_{\leq
-1}GD(M)$ is by assumption obtained by a finite number of extensions
of objects of the form $N'[s]$ with $N'$ weakly Koszul and $s>0$.
The module $H^{i+1}(F(N')[s])$ is cogenerated in degree $-i-s-1$, so
$(H^{i+1}(F(N')[s]))_{-i}=0$ for all $i \in \Z$ when $s>0$. By
induction $(H^{i+1}(F(\tau_{\leq -1}GD(M))))_{-i}=0$ for all $i \in
\Z$. But $H^i F(N)$ is cogenerated in degree $-i$, so $H^i F(N)
\simeq H^{i+1}(F(\tau_{\leq -1}GD(M)))=0$ when $i>0$.

So $H^i F(N) \neq 0$ only when $i=0$. This means that the weakly
Koszul module $N$ is generated in degree $0$ and is therefore
Koszul. So $\check E(M)$ has a direct summand $N=H^0 GD(M)$ which is
Koszul.
\end{proof}

Surprisingly, this is not the beginning of an inductive procedure.
An example in the next section (Example \ref{weak}) shows that the
other direct summands of $\check E(M)$ are not necessarily Koszul or
shifts of Koszul modules. In other words $\check E(M)$ can have
indecomposable direct summands which are generated in multiple
degrees.

\section{computation of the object $GD(M)$}\label{diff}

So far we have found a formula for the homology of the object
$GD(M)$ when $M$ is a weakly Koszul module, but we have not
described the object $GD(M)$ itself. Based on the description of the
functor $G$ in \cite{Bei}, we give a method for computing $GD(M)$
for any finitely generated $\L$-module when $\L$ is given as a path
algebra (that is an algebra given as a quiver with relations). We
refer to \cite{Dua} for more details on a construction that is
essentially the same as ours, but there it is done in an abelian
setting.

Suppose $\L$ is a Koszul algebra given as the path algebra of a
quiver with relations. Then $\L_0$ is the $k$-linear span of the
vertices, while $\L_1$ is the $k$-linear span of the arrows. It can
be shown that the relations are \emph{quadratic}, that is they are
given by a $\L_0$-sub-bimodule $R$ of $\L_1 \otimes_{\L_0} \L_1$.
The quiver of $\L^{\op}$ is the opposite quiver to the one of $\L$.
The relations for $\L^{\op}$ are similarly given by a
$\L_0^{\op}$-sub-bimodule $\check R$ of $\L_1^{\op}
\otimes_{\L_0^{\op}} \L_1^{\op}$. The Koszul dual $\G=\bigoplus_{i
\geq 0}\Ext_\L^i(\L_0,\L_0)$ is isomorphic to the path algebra with
the same quiver as $\L$ (so $\G_0 \simeq \L_0$, $\G_1 \simeq \L_1$
and $\G_1 \otimes_{\G_0} \G_1 \simeq \L_1 \otimes_{\L_0} \L_1$) and
\emph{orthogonal} relations $R^{\perp}$, that is the kernel of the
$\L_0$-bimodule morphism $\Hom_{\L_0^{\op}}(\L_1^{\op}
\otimes_{\L_0^{\op}} \L_1^{\op},\L_0^{\op}) \pil
\Hom_{\L_0^{\op}}(\check R,\L_0^{\op})$ (see \cite{Bei}). For path
algebras, the orthogonal relations can be found with the help of the
bilinear form used in \cite[2.2]{Gre}.

Now let $Q$ denote the quiver of $\L$, let $\check Q$ denote the
quiver of $\L^{\op}$ and let $Q^{\ast}$ denote the quiver of $\G$.
Here $Q^{\ast}=Q$ but it is still helpful to keep separate notation.
If the vertices of $Q$ are indexed by $\{ 1, \ldots, t \}$, let $\{
\check 1, \ldots, \check t \}$ denote the corresponding vertices of
$\check Q$ and let $\{ 1^{\ast}, \ldots, t^{\ast} \}$ denote the
corresponding vertices of $Q^{\ast}$. If $a \extto \alpha b$ is an
arrow in $Q$, then let $\check a$ denote the corresponding arrow
$\check a \xleftarrow {\check \alpha} \check b$ in $\check Q$ let
$\alpha^{\ast}$ denote the corresponding arrow $a^{\ast} \extto
{\alpha^{\ast}} b^{\ast}$ in $Q^{\ast}$. If $\check a$ is a vertex
of $\check Q$, let $S_{\check a}$ denote the corresponding simple
$\L^{\op}$-module (and $\L^{\op}_0$-module). Let $P_{a^{\ast}}$
denote the projective $\G$-module corresponding to the vertex
$a^{\ast}$ of $Q^{\ast}$.

Now let $M$ be a finitely generated $\L$-module. For the moment we
do not assume that $M$ is weakly Koszul. Then $DM$ is a finitely
cogenerated $\L^{\op}$-module. We now try to find a complex
representing $GD(M)$, so we need to know the result of applying $G$
to a finitely cogenerated module.

The terms of $GD(M)$ we find from the graded parts of $DM$ (or just
as easily directly from the graded parts of $M$). If $(DM)_i \simeq
(S_{\check 1})^{n_1} \bigoplus \ldots \bigoplus (S_{\check
t})^{n_t}$ is a decomposition of $(DM)_i$ into simple
$\L^{\op}_0$-modules, then we put $GD(M)^i=[(P_{1^\ast})^{n_1}
\bigoplus \ldots \bigoplus (P_{t^\ast})^{n_t}] \gsh {-i}$.

To describe the differential is slightly more complicated. Suppose
$\check a$ and $\check b$ are two (not necessarily distinct)
vertices in $\check Q$ with $r$ arrows going from $\check b$ to
$\check a$. Denote the arrows by $\check \alpha_1, \ldots, \check
\alpha_r$.
$$\xymatrix{\check a & \ar@<-3.5ex>[l]_{\check \alpha_1}
\ar@<-1ex>[l]_{\check \alpha_2} \ar@{}[l]|{\vdots}
\ar@<2.5ex>[l]^{\check \alpha_r} \check b}$$ (There might also be
arrows in the other direction.) Suppose $x$ is an element in the
summand $(S_{\check b})^{n_{b,i}}$ of $(DM)_i$ and let $\check a
\xleftarrow {\check \alpha_j} \check b$ be an arrow. Then $\check
\alpha_j \cdot x$ is an element in the summand $(S_{\check
a})^{n_{a,i+1}}$ of $(DM)_{i+1}$. Choose $k$-bases for $(S_{\check
b})^{n_{b,i}}$ and $(S_{\check a})^{n_{a,i+1}}$ and let $A_{(\check
\alpha_j,i)}$ be the $(n_{a,i+1} \times n_{b,i})$-matrix with
entries in $k$ which represents the map $(S_{\check b})^{n_{b,i}}
\pil (S_{\check a})^{n_{a,i+1}}$ induced by $\check \alpha_j$.

For each arrow $a^{\ast} \extto {\alpha_j^{\ast}} b^{\ast}$ and each
$i \in \Z$ we have a map $(\bar {\alpha_j^{\ast}})_i \colon
P_{b^\ast} \gsh {-i} \pil P_{a^\ast} \gsh {-i-1}$ which we can view
as multiplication with $\alpha_j^{\ast}$ from the right. The part of
$d^i \colon GD(M)^i \pil GD(M)^{i+1}$ which maps $(P_{b^\ast} \gsh
{-i})^{n_{b,i}}$ to $(P_{a^\ast} \gsh {-i-1})^{n_{a,i+1}}$ is given
by $$\sum_{j=1}^r (A_{(\check \alpha_j,i)} \times (\bar
{\alpha_j^{\ast}})_i),$$ where $\times$ means that each entry in the
matrix $A_{(\check \alpha_j,i)}$ is to be multiplied by $(\bar
{\alpha_j^{\ast}})_i$. We illustrate with an example.

\begin{example}\label{noeth}
Let $\L$ be the path algebra of the quiver
$$\xymatrix{1 \ar[r]^{\alpha} & 2  \ar@(ur,dr)^{\beta}}.$$ This
is a Koszul algebra.

Then $\L^{\op}$ is the path algebra of the quiver
$$\xymatrix{\check 1 & 2 \ar[l]_{\check \alpha} \ar@(dr,ur)_{\check \beta}},$$
and $\G \simeq \bigoplus_{i \geq 0}\Ext^i_\L(\L_0,\L_0)$ is
isomorphic to the path algebra given by the quiver
$$\xymatrix{1^{\ast} \ar[r]^{\alpha^{\ast}} & 2^{\ast}  \ar@(ur,dr)^{\beta^{\ast}}}$$
and with relations $\beta^{\ast} \alpha^{\ast}=0$ and
$(\beta^{\ast})^2=0$.

Let $M$ be the following infinite dimensional $\L$-module generated
in degrees $0$, $1$ and $2$
$$\xymatrix @=1.5ex{S_1 \ar@{-}[dr] \\ S_1
\ar@{-}[dr] & S_2 \ar@{-}[d] \\ S_1 \ar@{-}[dr] & S_2 \ar@{-}[d]
\\& S_2 \ar@{-}[d]\\
& { } \ar@{.}[d] \\& {} }$$ Here $\tmfgr (M) \simeq P_1 \oplus S_1
\gsh 1 \oplus S_1 \gsh 2$. Since $\mfgr (M) \simeq P_1 \oplus
(S_1)^2$ is a Koszul module, the module $M$ is weakly Koszul.

Then $DM$ is the finitely cogenerated $\L^{\op}$-module
$$\xymatrix @=1.5ex{& {} \ar@{.}[d] \\ & S_{\check 2} \ar@{-}[d] \\ &
S_{\check 2} \ar@{-}[dl] \ar@{-}[d]
\\
S_{\check 1} & S_{\check 2} \ar@{-}[dl] \ar@{-}[d] \\
S_{\check 1} & S_{\check 2} \ar@{-}[dl] \\ S_{\check 1}}$$

From this we can read off the object $GD(M)$: $$\ldots \pil
P_{2^{\ast}} \gsh 4 \extto {\bar \beta^{\ast}} P_{2^{\ast}} \gsh 3
\extto {\begin{pmatrix}  \bar \alpha^{\ast} \\
\bar \beta^{\ast}
\end{pmatrix}} (P_{1^{\ast}} \oplus
P_{2^{\ast}}) \gsh 2 \extto {\begin{pmatrix} 0 & \bar \alpha^{\ast} \\
0 & \bar \beta^{\ast}
\end{pmatrix}} (P_{1^{\ast}} \oplus P_{2^{\ast}}) \gsh 1
\extto {\begin{pmatrix} 0 & \bar \alpha^{\ast}
\end{pmatrix}} P_{1^{\ast}}
\pil 0$$ Here \begin{align*} H^0 GD(M) &\simeq S_{1^{\ast}}  \simeq
\check E(P_1) \simeq H^0 GD(\tmfgr (M)),\\  H^{-1} GD(M) &\simeq
P_{1^{\ast}} \gsh 1  \simeq \check E(S_1) \gsh 1 \simeq H^{-1}
GD(\tmfgr (M)),\\ H^{-2} GD(M) &\simeq P_{1^{\ast}} \gsh 2  \simeq
\check E(S_1) \gsh 2 \simeq H^{-2} GD(\tmfgr (M))\\ \intertext{and}
H^n GD(M) &\simeq 0 \simeq H^n GD(\tmfgr (M))
\end{align*} for $n \notin \{-2,-1,0\}$ as expected.
\end{example}

In the next example we compute $GD(M)$ for a module $M$ which is not
weakly Koszul but turns out to be quasi-Koszul.

\begin{example}\label{quasi}
Let $\L$ be the path algebra of the quiver $$\xymatrix{&& 5
\ar[d]_\delta \ar@{--}@/_1pc/ [dr] & \\1 \ar[r]_\alpha & 2
\ar[r]_\beta & 3 \ar[r]_\gamma  & 4 }$$ and relation $\gamma
\delta=0$. This algebra is Koszul.

Then $\G \simeq \bigoplus_{i \geq 0}\Ext^i_\L(\L_0,\L_0)$ is
isomorphic to the path algebra of the quiver $$\xymatrix{&& 5^{\ast}
\ar[d]_{\delta^{\ast}} & \\1^{\ast} \ar[r]^{\alpha^{\ast}}
\ar@{--}@/_1pc/ [rr]& 2^{\ast} \ar[r]^{\beta^{\ast}} \ar@{--}@/_1pc/
[rr]& 3^{\ast} \ar[r]^{\gamma^{\ast}} & 4^{\ast} }$$ \newline with
relations $\beta^{\ast} \alpha^{\ast}=0$ and $\gamma^{\ast}
\beta^{\ast}=0$.

Let $M$ be the following module generated in degrees $0$ and $1$.
$$
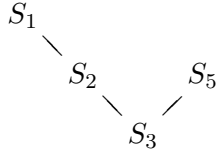
\xymatrix @R=1.5ex @C=1ex{S_1 \ar@{-}[dr]&&&
\\& S_2 \ar@{-}[dr] && S_5\ar@{-}[dl]\\&& S_3 & }$$
Since $\mfgr (M) \simeq I_3 \oplus S_5$ is not Koszul, this is not a
weakly Koszul module.

Then $GD(M)$ is represented by the complex
$$0 \pil P_{3^{\ast}} \gsh 2 \extto {\begin{pmatrix}  \bar \beta^{\ast} \\
\bar \delta^{\ast}
\end{pmatrix}} (P_{2^{\ast}} \oplus P_{5^{\ast}}) \gsh 1 \extto {\begin{pmatrix} \bar
\alpha^{\ast} & 0
\end{pmatrix}} P_{1^{\ast}} \pil 0$$

Here $H^n GD(M)=0$ for $n \notin \{-1,0\}$, $H^0 GD(M) \simeq
S_{1^{\ast}}$ and $H^{-1} GD(M)$ is the module
$$
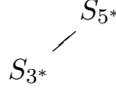
\xymatrix @R=1.5ex @C=1ex{& S_{5^{\ast}} \ar@{-}[dl]\\ S_{3^{\ast}} &
}$$ generated in degree $1$. In this example, although $H^i GD(M)$
is generated in degree $-i$ for all $i \in \Z$, the module $H^{-1}
GD(M)$ does not have a linear resolution. Therefore the conditions
in Theorem \ref{char} are not satisfied. Since $\check E (M) \simeq
\bigoplus_{i \in \Z} (H^i GD(M)) \gsh i$ is generated in degree $0$,
we have that $M$ is quasi-Koszul.
\end{example}

In Theorem \ref{eweak} we have shown that if $\check E(M)$ is weakly
Koszul, then $\check E(M)$ has a Koszul direct summand. The
following example shows that $\check E(M)$ can have other
indecomposable direct summands which are generated in multiple
degrees.

\begin{example}\label{weak}
Let $\L$ be the path algebra of the quiver $$
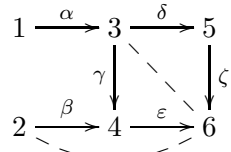
\xymatrix{1
\ar[r]^{\alpha} & 3 \ar[r]^{\delta} \ar[d]_{\gamma} \ar@{--}[dr] & 5
\ar[d]^{\zeta}
\\2 \ar[r]^{\beta} \ar@{--}@/_1pc/ [rr] & 4 \ar[r]^{\varepsilon} & 6}$$ and relations
$\varepsilon \beta=0$ and $\varepsilon \gamma - \zeta \delta=0$.
This algebra is Koszul.

Then $\G \simeq \bigoplus_{i \geq 0}\Ext^i_\L(\L_0,\L_0)$ is
isomorphic to the path algebra of the quiver $$
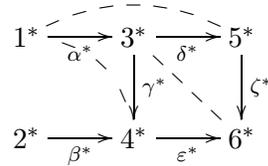
\xymatrix{1^{\ast}
\ar[r]_{\alpha^{\ast}} \ar@{--}@/^1pc/ [dr] \ar@{--}@/^1pc/ [rr] &
3^{\ast} \ar[r]_{\delta^{\ast}} \ar[d]^{\gamma^{\ast}} \ar@{--}[dr]
 & 5^{\ast} \ar[d]^{\zeta^{\ast}}
\\2^{\ast} \ar[r]_{\beta^{\ast}} & 4^{\ast} \ar[r]_{\varepsilon^{\ast}} & 6^{\ast}}$$
with relations $\gamma^{\ast} \alpha^{\ast}=0$, $\delta^{\ast}
\alpha^{\ast}=0$ and $\varepsilon^{\ast} \gamma^{\ast} +
\zeta^{\ast} \delta^{\ast}=0$.

Let $M$ be the following $\L$-module generated in degrees $0$ and
$1$.
$$
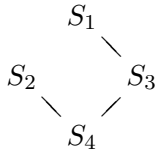
\xymatrix @R=1.5ex @C=1ex{& S_1 \ar@{-}[dr] &
\\S_2 \ar@{-}[dr] && S_3\ar@{-}[dl]\\& S_4 & }$$

Then $GD(M)$ is represented by the complex
$$0 \pil P_{4^{\ast}} \gsh 2 \extto {\begin{pmatrix}  \bar \beta^{\ast} \\
\bar \gamma^{\ast}
\end{pmatrix}} (P_{2^{\ast}} \oplus P_{3^{\ast}}) \gsh 1 \extto {\begin{pmatrix} 0 & \bar
\alpha^{\ast}
\end{pmatrix}} P_{1^{\ast}} \pil 0$$

Here $H^n GD(M)=0$ for $n \notin \{-1,0\}$, $H^0 GD(M) \simeq
S_{1^{\ast}}$ and $H^{-1} GD(M)$ is the $\G$-module
$$\xymatrix @R=1.5ex @C=1ex{S_{2^{\ast}} \ar@{-}[dr]&&&
\\& S_{4^{\ast}} \ar@{-}[dr] && S_{5^{\ast}} \ar@{-}[dl]\\&& S_{6^{\ast}} & }$$
generated in degrees $1$ and $2$. Call this module $L$. Since $\mfgr
(L) \simeq P_{2^{\ast}} \oplus S_{5^{\ast}}$ is Koszul, the module
$L$ is weakly Koszul. From Lemma \ref{hatem} we get $\check E(M)
\simeq L \gsh {-1} \oplus S_{1^{\ast}}$. So $\check E(M)$ is weakly
Koszul and has an indecomposable summand which is generated in
multiple degrees.
\end{example}

\section{$\ainf$-modules}\label{ainf}

In this section we discuss an alternative way of viewing objects in
$\Dgr \G$, namely as $\ainf$-modules. Instead of thinking of objects
as complexes, we think of them as homology groups with some
additional structure. If we fix homology groups satisfying the
conditions of Theorem \ref{char}, then each possible $\ainf$-module
structure on that homology gives an object in $\X$. From Corollary
\ref{sist} we know that if two objects in $\X$ share the same
homology, then the two corresponding weakly Koszul $\L$-modules have
isomorphic adjusted associated graded modules. So if we classify all
objects in $\X$ with a certain homology, then via Koszul duality we
classify all weakly Koszul $\L$-modules with a certain adjusted
associated graded module.

Let $\G=\bigoplus_{i \geq 0}\G_i$ be a graded algebra. We consider
$\G$ as an $\ainf$-algebra concentrated in degree $0$. The ordinary
grading of $\G$ gives an additional structure which is also
inherited by our modules. So what we really are considering are
\emph{graded} $\ainf$-modules. This extra grading can be introduced
more formally by considering $\ainf$-algebras (and their
$\ainf$-modules) over the monoidal base category of graded vector
spaces, but we will not take this approach here. For definitions we
follow \cite{Doc} and \cite{Kel}.

For us a \emph{(graded) $\ainf$-module} over $\G$ is a bigraded
space $$N=\bigoplus_{(i,j) \in \Z \times \Z}N_j^i$$ with maps
$$m_n \colon \G^{\otimes n-1} \otimes N \pil N,\textit{ } n \geq 1$$
of bidegree $(n-2,0)$ satisfying the rules
\begin{align*}
m_1 m_1 &=0,\\
m_1 m_2 &=m_2(\mathtt 1 \otimes m_1)\\
\intertext{and for $n \geq 3$}
\sum_{i=1}^n(-1)^{i(n-1)}m_{n-i+1}(\mathtt 1^{\otimes n-i} \otimes
m_i) &=\sum_{j=1}^{n-2}(-1)^{j-1} m_{n-1}(\mathtt 1^{\otimes n-j-2}
\otimes m \otimes \mathtt 1^{\otimes j})
\end{align*}
where $\mathtt 1$ is the identity map and $m$ is the multiplication
of $\G$. Some terms are omitted from the usual definition since
there are no higher multiplications in $\G$.

We only consider \emph{strictly unital} modules, that is modules $N$
such that for all $a \in N$, we have $m_2(1,a)=a$ and
$m_n(\gamma_1,\ldots,\gamma_{n-1},a)=0$ if $n \geq 3$ and $1 \in
\{\gamma_1,\ldots,\gamma_{n-1} \}$.

Note that if we let $N^i=\bigoplus_{j \in \Z} N_j^i$, then the $N^i$
together with $m_1$ form a complex $(N^{\bullet},m_1)$ of graded
$k$-modules.

Two special cases are important. The first is when $m_n=0$ for $n
\geq 3$, the other is when $m_1=0$. In the second case each $N^i$ is
a graded $\G$-module, not only a graded $k$-module. In the first
case $N$ is essentially a complex of graded $\G$-modules, and we
view complexes of graded $\G$-modules in this way.

A \emph{morphism} $f \colon L \pil N$ between two $\ainf$-modules
$L$ and $N$ is given by a family of maps
$$f_n \colon \G^{\otimes n-1} \otimes L \pil N, \textit{ } n \geq
1$$ of bidegree $(n-1,0)$ satisfying the rules
\begin{align*}
f_1 m_1 &=m_1 f_1,\\
f_1 m_2-f_2(\mathtt 1 \otimes m_1) &=m_2(\mathtt 1 \otimes f_1)+m_1
f_2\\
\intertext{and for $n \geq 3$}
\sum_{i=1}^n(-1)^{i(n-1)}f_{n-i+1}(\mathtt 1^{\otimes n-i} \otimes
m_i) &+\sum_{j=1}^{n-2}(-1)^j f_{n-1}(\mathtt 1^{\otimes n-j-2}
\otimes m \otimes \mathtt 1^{\otimes j})\\ &=\sum_{r=1}^n
(-1)^{(r+1)n} m_{n-r+1}(\mathtt 1^{\otimes n-r} \otimes f_r).
\end{align*}
Note in particular that $f_1$ is a chain map $f_1 \colon
(L^{\bullet},m_1) \pil (N^{\bullet},m_1)$ between complexes of
graded $k$-modules.

We only consider \emph{strictly unital} morphisms, that is morphisms
$f$ such that $f_n(\lambda_1, \ldots,\lambda_{n-1},a)=0$ whenever $n
\geq 2$ and $1 \in \{\lambda_1,\ldots,\lambda_{n-1} \}$.

The \emph{identity morphism} $f \colon N \pil N$ is given by
$f_1=\mathtt 1$ and $f_i=0$ for all $i>0$. The composition $fg
\colon N \pil M$ of two morphisms $f \colon L \pil M$ and $g \colon
N \pil L$ is given by the rule
$$(fg)_n=\sum_{i=1}^n f_{n-i+1}(\mathtt 1^{\otimes n-i} \otimes
g_i).$$

We say that $f$ is a \emph{quasi-isomorphism} if $f_1$ is a
quasi-isomorphism. With the definitions we have made, the
quasi-isomorphism classes we get do not differ significantly from
the ones we have for complexes of graded $\G$-modules. Each
quasi-isomorphism class of $\ainf$-modules over $\G$ corresponds to
(has as a subclass) exactly one quasi-isomorphism class of complexes
of graded $\G$-modules.

An important theorem \cite[3.3.1.7]{Doc} states that for any
strictly unital $\ainf$-module $N$, there is a $\ainf$-module
structure on $H^{\ast}N$ with $m_1=0$ and a strictly unital
quasi-isomorphism between $H^{\ast}N$ with this structure and $N$.
If $L$ and $N$ are two modules, both with $m_1=0$, then each
quasi-isomorphism between them is an isomorphism.

Therefore if we want to describe an object in $\Dgr \G$ (or strictly
speaking an equivalent category), it suffices to specify its
homology and an $\ainf$-module structure in its quasi-isomorphism
class with $m_1=0$. If we want to classify all objects (up to
isomorphism) in $\X$ with a certain homology, it suffices to
classify all isomorphism classes of possible $\ainf$-module
structures with $m_1=0$ on this homology.

\begin{example}
Let $\L$ and $\G$ be the same as in Example \ref{noeth}. Let $t \geq
1$ be a number and let $M_{(t)}$ denote the following infinite
dimensional $\L$-module generated in degrees $0$ and $t$.
$$\xymatrix @=1.5ex{S_1 \ar@{-}[dr] \\ & S_2 \ar@{-}[d] \\
&{ } \ar@{.}[d]\\ &{ } \ar@{-}[d]\\ S_1 \ar@{-}[dr] & S_2 \ar@{-}[d]
\\& S_2 \ar@{-}[d]\\
& { } \ar@{.}[d] \\& {} }$$

In this example we find an $\ainf$-module structure (with $m_1=0$)
for $GD(M_{(t)})$.

Since $\tmfgr(M_{(t)}) \simeq P_1 \oplus S_1 \gsh t$ is a Koszul
module, the module $M_{(t)}$ is weakly Koszul. Therefore by
Corollary \ref{form} we have
\begin{align*} H^0 GD(M_{(t)}) &\simeq  \check E(P_1) \simeq S_{1^{\ast}},\\
H^{-t} GD(M_{(t)}) &\simeq \check E(S_1) \gsh t \simeq P_{1^{\ast}}
\gsh t \intertext{and} H^n GD(M_{(t)}) &\simeq 0
\end{align*} for $n \notin \{-t,0\}$.

We now look for possible $\ainf$-structures (with $m_1=0$) on this
homology. Fix a basis vector $v$ for $S_{1^{\ast}}$ and a basis
vector $w$ for the socle of $P_{1^{\ast}} \gsh t$. Due to degree
considerations (and the remark following this example), the only
possibly non-zero higher structure is that
$m_{t+2}(\beta,\beta,\ldots,\beta,\alpha,v)=xw$ for some $x \in k$.
All values of $x$ give permissible $\ainf$-structures.

Again due to degree considerations, all quasi-isomorphisms between
such structures must have $f_i=0$ for $i \geq 2$. It is possible and
easy to construct quasi-isomorphisms using only $f_1$ between
structures with $x \neq 0$. We choose $x=1$ as a representative for
this orbit and denote by $N_{(t)}$ the corresponding object in $\Dgr
\G$. The remaining case is $x=0$ and corresponds to the object
$$P_{1^{\ast}} \gsh t [t] \oplus S_{1^{\ast}} \simeq GD(P_1 \oplus S_1 \gsh t).$$
Since $M_{(t)} \not \simeq P_1 \oplus S_1 \gsh t$, we must have
$GD(M_{(t)}) \simeq N_{(t)}$.

Therefore $GD(M_{(t)})$ can described as the homology
\begin{align*} H^0 GD(M_{(t)}) &\simeq S_{1^{\ast}},\\
H^{-t} GD(M_{(t)}) &\simeq P_{1^{\ast}} \gsh t \end{align*} with
additional $\ainf$-structure
$$m_{t+2}(\beta,\beta,\ldots,\beta,\alpha,v)=w.$$

Since $t$ can be chosen arbitrarily large and $m_{t+2} \neq 0$, this
example shows that arbitrarily high module structures are needed to
describe all objects in $\Dgr \G$ in this way.

The weakly Koszul $\L$-modules with adjusted associated graded
module $P_1 \oplus S_1 \gsh t$ we have found to be $P_1 \oplus S_1
\gsh t$ and $M_{(t)}$.
\end{example}

\begin{remark}
We are assuming that the higher module structure maps $m_n$ not only
respect the ordinary grading, but also the grading given by the
quiver. This can also be justified by a change of monoidal base
category.
\end{remark}

\begin{example}
In this example we show how to find all objects $N$ in $\Dgr \G$
with the same homology as $GD(M)$ in Example \ref{noeth}. Each
isomorphism class corresponds to a weakly Koszul module $\tilde M$
over $\L$ with $\tmfgr (\tilde M) \simeq \tmfgr (M)$. As a result of
our computation we also find an $\ainf$-module structure for
$GD(M)$.

The given homology is $H^0 N \simeq S_{1^{\ast}}$, $H^{-1} N \simeq
P_{1^{\ast}} \gsh 1$, $H^{-2} N \simeq P_{1^{\ast}} \gsh 2$ and $H^n
N \simeq 0$ for $n \notin \{-2,-1,0\}$. In total this homology is
$5$-dimensional. We fix basis vectors
\begin{align*}v_1
&\in (H^0 N)_0,\\ v_2 &\in (H^{-1} N)_1,\\ v_3 &\in (H^{-1} N)_2,\\
v_4 &\in (H^{-2} N)_2,\\ v_5 &\in (H^{-2} N)_3.
\end{align*}

The possible higher products are $m_3(\beta,\alpha,v_1)=x v_3$,
$m_3(\beta,\alpha,v_2)=y v_5$ and $m_4(\beta,\beta,\alpha,v_1)=z
v_5$ where $x$, $y$ and $z$ are elements in $k$. All triples
$(x,y,z)$ give permissible $\ainf$-structures, so we have a
$3$-dimensional representation space. We now want to find the
isomorphism classes. The possible quasi-isomorphisms are given by
$f_1(v_i)=q_i v_i$, $1 \leq i \leq 5$ and  $f_2(\beta,v_3)=\mu v_5$,
where $q_i,\mu \in k$ and $q_i \neq 0$, $1 \leq i \leq 5$. It
follows from the formulas that $q_2=q_3$ and $q_4=q_5$. If we let
$a=q_2/q_1$, $b=q_4/q_2$ and $\rho =\mu/q_1$, then possible
quasi-isomorphisms between triples $(x,y,z)$ and $(x',y',z')$ are
given by the formulas
\begin{align*}
x' &=ax,\\
y' &=by+\rho,\\
z' &=abz-a\rho x
\end{align*}

This divides the representation space into $4$ orbits, namely
\begin{align*}
\mathcal{O}_1 &=\{(0,\rho,0) \mid \rho \in k \},\\
\mathcal{O}_2 &=\{(a,\rho,-a \rho) \mid a,\rho \in k; a \neq 0\},\\
\mathcal{O}_3 &=\{(0,\rho,c) \mid c, \rho \in k; c \neq 0\},\\
\mathcal{O}_4 &=\{(a,\rho,c-a \rho) \mid a,c, \rho \in k; a, c \neq
0\}
\end{align*}

We choose representatives $\mathbf{x}_1=(0,0,0)$,
$\mathbf{x}_2=(1,0,0)$, $\mathbf{x}_3=(0,0,1)$ and
$\mathbf{x}_4=(1,0,1)$ respectively. The object in $\Dgr \G$
corresponding to $\mathbf x_1$ is \begin{align*}S_{1^{\ast}} \oplus
P_{1^{\ast}} \gsh 1 [1] \oplus P_{1^{\ast}} \gsh 2 [2]&\simeq GD(P_1
\oplus S_1 \gsh 1 \oplus S_1 \gsh 2).\\ \intertext{From the previous
example we recognise the object corresponding to $\mathbf x_2$ as}
N_{(2)} \oplus P_{1^{\ast}} \gsh 2 [2]&\simeq GD(M_{(2)} \oplus S_1
\gsh 2), \intertext{and the object correponding to $\mathbf x_3$ as}
N_{(3)} \oplus P_{1^{\ast}} \gsh 1 [1]&\simeq GD(M_{(3)} \oplus S_1
\gsh 1).\end{align*} Since $GD(M)$ is not isomorphic to any of
these, it must correspond to $\mathbf x_4$.

So $GD(M)$ can be described as the homology
\begin{align*} H^0 GD(M) &\simeq S_{1^{\ast}},\\
H^{-1} GD(M) &\simeq P_{1^{\ast}} \gsh 1,\\
H^{-2} GD(M) &\simeq P_{1^{\ast}} \gsh 2 \end{align*} with
additional $\ainf$-structure
\begin{align*} m_3(\beta,\alpha,v_1)=v_3,\\
m_4(\beta,\beta,\alpha,v_1)=v_5.\end{align*}

The weakly Koszul $\L$-modules with adjusted associated graded
module $P_1 \oplus S_1 \gsh 1 \oplus S_1 \gsh 2$ we have found to be
$P_1 \oplus S_1 \gsh 1 \oplus S_1 \gsh 2$, $M_{(2)} \oplus S_1 \gsh
2$, $M_{(3)} \oplus S_1 \gsh 1$ and $M$.

\end{example}

\end{document}